\documentclass[a4paper,12pt]{article}

\usepackage{amsmath,amstext,amscd,amsthm,amsfonts}
\usepackage{amssymb}
\newtheorem*{theo}{Proposition}
\newtheorem{theorem}{Theorem}[section]
\newtheorem{lemma}[theorem]{Lemma}
\newtheorem{proposition}[theorem]{Proposition}

\newtheorem*{conjecture}{Conjecture}
\theoremstyle{remark}
\newtheorem{rk}[theorem]{Remark}

\def\Ad{\mathop{\rm Ad}\nolimits}
\def\ad{\mathop{\rm ad}\nolimits}

\def\Mat{\mathop{\rm Mat}\nolimits}
\def\Mod{\mathop{\rm mod}\nolimits}

\def\Lie{\mathop{\rm Lie}\nolimits}
\def\GL{\mathop{\rm GL}\nolimits}

\def\SO{\mathop{\rm SO}\nolimits}

\def\Sp{\mathop{\rm Sp}\nolimits}

\def\dim{\mathop{\rm dim}\nolimits}

\def\End{\mathop{\rm End}\nolimits}

\def\rank{\mathop{\rm rk}\nolimits}
\def\min{\mathop{\rm min}\nolimits}

\def\rank{\mathop{\rm rk}\nolimits}

\def\diag{\mathop{\rm diag}\nolimits}

\def\triv{\mathop{\rm triv}\nolimits}

\title{Varieties of Modules for ${\mathbb Z}/2{\mathbb Z}\times{\mathbb Z}/2{\mathbb Z}$}

\author{Paul Levy \\
Ecole Polytechnique F\'ed\'erale de Lausanne \\
EPFL SB IGAT \\
B\^atiment BCH \\
1015 Lausanne \\
Switzerland \\
paul.levy@epfl.ch}

\begin{document}

\maketitle{}{}

\begin{abstract}
Let $k$ be an algebraically closed field of characteristic 2.
We prove that the restricted nilpotent commuting variety ${\mathcal C}$, that is the set of pairs of $(n\times n)$-matrices $(A,B)$ such that $A^2=B^2=[A,B]=0$, is equidimensional.
${\mathcal C}$ can be identified with the `variety of $n$-dimensional modules' for ${\mathbb Z}/2{\mathbb Z}\times{\mathbb Z}/2{\mathbb Z}$, or equivalently, for $k[X,Y]/(X^2,Y^2)$.
On the other hand, we provide an example showing that the restricted nilpotent commuting variety is not equidimensional for fields of characteristic $>2$.
We also prove that if $e^2=0$ then the set of elements of the centralizer of $e$ whose square is zero is equidimensional.
Finally, we express each irreducible component of ${\mathcal C}$ as a direct sum of indecomposable components of varieties of ${\mathbb Z}/{2{\mathbb Z}}\times{\mathbb Z}/2{\mathbb Z}$-modules.
\end{abstract}

\section{Introduction}

Let $G=\GL(n,k)$, where $k$ is an algebraically closed field of characteristic $p>0$, and let ${\mathfrak g}$ be the Lie algebra of $G$.
We denote the $p$-th power of matrices on ${\mathfrak g}$ by $x\mapsto x^{[p]}$, and its iteration $m$ times by $x\mapsto x^{[p^m]}$.
(This is the standard notation in the theory of restricted Lie algebras.)
Clearly $x$ is nilpotent if and only if $x^{[p^N]}=0$ for $N\gg 0$.
Denote by ${\cal N}$ the set of nilpotent elements of ${\mathfrak g}$ and by ${\cal N}_1$ the subset of elements satisfying $x^{[p]}=0$ (the {\it restricted nullcone}).
It was proved in \cite{npv} that ${\cal N}_1$ is irreducible.
(An explicit description was given in \cite{clnp}.)
In \cite{premnil}, Premet proved that the nilpotent commuting variety ${\cal C}^{nil}({\mathfrak g}):=\{ (x,y)\in{\mathfrak g}\times{\mathfrak g}\,|\, x,y\in{\cal N},[x,y]=0\}$ is irreducible and of dimension $(n^2-1)$.
More specifically, ${\cal C}^{nil}({\mathfrak g})=\overline{G\cdot(e,{\mathfrak u})}$ where $e$ is a regular nilpotent element of ${\mathfrak g}$ and ${\mathfrak u}=ke\oplus ke^2\oplus\ldots\oplus ke^{n-1}$.
(Here and in what follows we use the dot to denote the action of $G$ by conjugation on ${\mathfrak g}$ or the induced diagonal action on ${\mathfrak g}\times{\mathfrak g}$, and the notation $\overline{V}$ for the Zariski closure of a subset $V$ of an arbitrary affine vector space, where the context is clear.)
In fact, Premet proved that the nilpotent commuting variety of $\Lie(G)$ is equidimensional for any reductive group $G$ over an algebraically closed field of good characteristic.

The nilpotent commuting variety, or more accurately the restricted nilpotent commuting variety ${\cal C}^{nil}_1({\mathfrak g})=\{ (x,y)\in{\cal N}_1\times{\cal N}_1\,:\, [x,y]=0\}$ is related to the cohomology of $G$ by work of Suslin, Friedlander \& Bendel.
It was proved in \cite{sfb} that ${\cal C}^{nil}_1({\mathfrak g})$ is homeomorphic to the spectrum of the cohomology ring $\oplus_{i\geq 0} H^{2i}(G_2,k)$, where $G_2$ is the second Frobenius kernel of $G$.
More generally, the restricted nullcone ${\cal N}_1$ plays an important role in the representation theory of ${\mathfrak g}$ due to the theory of support varieties of (reduced enveloping algebras of) restricted Lie algebras (studied for the restricted enveloping algebra in \cite{fp1,fp2,jantzen} and for general reduced enveloping algebras in \cite{fp3}; see also \cite{premsupp,premcomp}).

Another perspective is that of varieties of modules (see for example \cite{cbs}): ${\cal C}_1^{nil}({\mathfrak g})$ can be identified with the variety of $n$-dimensional modules for the truncated polynomial ring $k[X,Y]/(X^p,Y^p)$.
There is an isomorphism $k[X,Y]/(X^p,Y^p)\rightarrow k\Gamma$, where $\Gamma={\mathbb Z}/p{\mathbb Z}\times{\mathbb Z}/p{\mathbb Z}$.
Specifically, if $\sigma$ (resp. $\tau$) is a generator for the first (resp. second) copy of ${\mathbb Z}/p{\mathbb Z}$ in $\Gamma$, then $X+1\mapsto \sigma$ and $Y+1\mapsto\tau$.

This paper began as a preliminary investigation of the restricted nilpotent commuting variety in the simplest possible case: hence we assume $p=2$.
In Section 1 we show that projection onto the first coordinate maps any irreducible component of ${\cal C}^{nil}_1({\mathfrak g})$ onto ${\cal N}_1$ (that is, the components are `determined' by the dense orbit in ${\cal N}_1$).
Equidimensionality then follows by equidimensionality of ${\mathfrak z}_{\mathfrak g}(e)\cap{\cal N}_1$, which we prove for any $e\in{\cal N}_1$.

\begin{theo}
Let $k$ be of characteristic 2 and let ${\cal C}_1^{nil}$ be the restricted nilpotent commuting variety.

(a) If $n=2m$, then ${\cal C}_1^{nil}$ has $([m/2]+1)$ irreducible components, each of dimension $3m^2$.

(b) If $n=2m+1$, then ${\cal C}_1^{nil}$ has $(m+1)$ irreducible components, each of dimension $3m(m+1)$.
\end{theo}

This result for ${\cal C}_1^{nil}$ might be expected to indicate that the restricted nilpotent commuting variety is equidimensional for general $p$.
However, we show that this is not the case (Remark \ref{notgeneral}).
On the other hand, we observe that ${\mathfrak z}_{\mathfrak g}(e)\cap{\cal N}_1$ is equidimensional for many choices of $G$ and $e$.
We conjecture that this is true for reductive $G$ in good characteristic.
We remark that the intersection ${\mathfrak z}_{\mathfrak g}(e)\cap{\cal N}_1$ can be identified with the support variety ${\cal V}_{{\mathfrak z}_{\mathfrak g}(e)}(k)$, where $k$ is the trivial ${\mathfrak z}_{\mathfrak g}(e)$-module.
In the final section we express each irreducible component of ${\cal C}_1^{nil}({\mathfrak g})$ as a direct sum of indecomposable components of modules.

Our method for obtaining the above results is a rather crude direct approach.
Such a strategy will clearly be inappropriate in general.

{\it Notation.}
We denote by $\Mat_{r\times s}$ the vector space of all $r\times s$ matrices over $k$.
If $x\in {\mathfrak g}$ then the centralizer of $x$ in ${\mathfrak g}$ (resp. $G$) will be denoted ${\mathfrak z}_{\mathfrak g}(e)$ (resp. $Z_G(e)$).
We will sometimes abuse notation and use ${\cal N}_1$ to refer to the set of $p$-nilpotent elements in an arbitrary Lie algebra.
This will cause no confusion.
We denote by $e_{ij}$ the matrix with 1 in the $(i,j)$-th position, and zeros everywhere else.
(The dimension of $e_{ij}$ will always be specified or clear from the context.)
Our convention is that all modules are left modules.
We denote by $[m/r]$ the integer part of the fraction $m/r$.

\section{Centralizers}

Let $G=\GL(n,k)$, let ${\mathfrak g}=\Lie(G)$ and let $e_0,e_1,\ldots ,e_m$ be a set of representatives for the orbits in ${\cal N}_1={\cal N}_1({\mathfrak g})$.
Clearly ${\cal C}^{nil}_1=\bigcup_{i=0}^m G\cdot(e_i,{\mathfrak z}_{\mathfrak g}(e_i)\cap{\cal N}_1)$.
In general the set ${\mathfrak z}_{\mathfrak g}(e_i)\cap{\cal N}_1$ is not irreducible.
For each $i$ let $V_i^{(1)},V_i^{(2)},\ldots ,V_i^{(r_i)}$ be the irreducible components of ${\mathfrak z}_{\mathfrak g}(e_i)\cap{\cal N}_1$.
The following Lemma is adapted from \cite[Prop. 2.1]{premnil}.
The argument works for arbitrary $G$ and $p$.
(The only requirement is that the number of orbits in ${\cal N}_1$ is finite.
This is well-known if $p$ is good (see \cite{rich}) but is true even if $p$ is bad \cite{holts}.)

\begin{lemma}\label{cpts}
Let $X$ be an irreducible component of ${\cal C}^{nil}_1$.
Then there is some $i$, $0\leq i\leq m$, and some $j$, $1\leq j\leq r_i$, such that $X=\overline{G\cdot(e_i,V_i^{(j)})}$.
Moreover, $V_i^{(j)}\subseteq \overline{G\cdot e_i}$.
\end{lemma}

\begin{proof}
Since there are finitely many of the sets $\overline{G\cdot(e_i,V_i^{(j)})}$ and they cover ${\cal C}_1^{nil}$, the first statement is obvious.
For the second statement, define an action of $\GL(2)$ on ${\mathfrak g}\times{\mathfrak g}$ by the morphism $\GL(2)\times({\mathfrak g}\times{\mathfrak g})\rightarrow{\mathfrak g}\times{\mathfrak g}$, $(\left(
\begin{array}{cc}
a & b \\
c & d
\end{array}
\right), (x,y))\mapsto(ax+by,cx+dy)$.
Clearly any element of $\GL(2)$ preserves ${\cal C}^{nil}_1$.
Hence $\GL(2)$ preserves each irreducible component of ${\cal C}^{nil}_1$.
In particular, $\tau(X)=X$, where $\tau:(x,y)\mapsto (y,x)$.
Suppose therefore that $X=\overline{G\cdot(e_i,V_i^{(j)})}$ is an irreducible component of ${\cal C}_1^{nil}$.
Let $\pi:{\mathfrak g}\times{\mathfrak g}\rightarrow{\mathfrak g}$ denote the first projection.
Then $\overline{\pi(X)}=\overline{G\cdot e_i}$.
But $X=\tau(X)$, hence $V_i^{(j)}\subseteq \overline{\pi(X)}$.
\end{proof}

Suppose from now on that $p=2$.
For each $i$, $0\leq i\leq m=[n/2]$, let $e_i=\left(
\begin{array}{ccc}
0 & 0 & I_i \\
0 & 0 & 0 \\
0 & 0 & 0
\end{array}
\right)\in{\mathfrak g}$, where $I_i$ is the $i\times i$ identity matrix.
Here the top left, top right, bottom left and bottom right submatrices are $i\times i$, the top middle and bottom middle submatrices are $i\times (n-2i)$, the centre left and centre right submatrices are $(n-2i)\times i$, and the central submatrix is $(n-2i)\times (n-2i)$.
Then $\{ e_0,e_1,\ldots ,e_m\}$ is a set of representatives for the conjugacy classes in ${\cal N}_1$.
It is easy to see, with the standard description of nilpotent orbits via partitions of $n$, that $e_i$ corresponds to the partition $2^i.1^{n-2i}$.
Moreover, we have the following inclusions: $\{ 0\}=\overline{G\cdot e_0}\subset\overline{G\cdot e_1}\subset\ldots\subset\overline{G\cdot e_M}={\cal N}_1$.
The condition $V_i^{(j)}\subseteq\overline{G\cdot e_i}$ is clearly equivalent to the inequality: $\rank(y)\leq i$ for all $y\in V_i^{(j)}$.

Fix $i$ until further notice and let $x$ be an element of the centralizer ${\mathfrak z}_{\mathfrak g}(e_i)$, which must have the form
$$\left(
\begin{array}{ccc}
A & B & C \\
0 & E & F \\
0 & 0 & A
\end{array}
\right)\;:\;\;A,C\in\Mat_{i\times i} \; ,E\in\Mat_{(n-2i)\times(n-2i)},B\in\Mat_{i\times (n-2i)},F\in\Mat_{(n-2i)\times i}.$$

The requirement $x\in{\cal N}_1$ is equivalent, with this notation, to the conditions $A^2=BF+[A,C]=0$, $E^2=0$, $AB=BE$ and $EF=FA$.
We will frequently use $A,B,C,E,F$ to refer to these submatrices of (an arbitrary) $x\in{\mathfrak z}_{\mathfrak g}(e_i)$ where the element $x$ is clear from the context.
Assume for the rest of this section that $n-2i\geq 2$, and let $V$ be an irreducible component of ${\mathfrak z}_{\mathfrak g}(e_i)\cap{\cal N}_1$.
We shall prove that $V\not\subseteq\overline{G\cdot e_i}$.
We begin with the following lemma.

\begin{lemma}\label{Enonzero}
Suppose there is some element $x\in V$ such that $E\neq 0$.
Then $V$ is not contained in $\overline{G\cdot e_i}$.
\end{lemma}

\begin{proof}
Let the one-dimensional torus $\lambda:k^\times\rightarrow G$, $t\mapsto\left(
\begin{array}{ccc}
tI_i & 0 & 0 \\
0 & I_{n-2i} & 0 \\
0 & 0 & t^{-1}I_i
\end{array}
\right)$ act on ${\mathfrak g}$ by conjugation.
Since $\lambda(t)e_i\lambda(t^{-1})=t^2 e_i$, $\lambda(k^\times)$ preserves ${\mathfrak z}_{\mathfrak g}(e_i)$, and therefore preserves each irreducible component of ${\mathfrak z}_{\mathfrak g}(e_i)\cap{\cal N}_1$.
Thus if $x\in V$ then $x_0=\lim_{t\rightarrow 0}(\Ad\lambda(t)(x))=\left(
\begin{array}{ccc}
A & 0 & 0 \\
0 & E & 0 \\
0 & 0 & A
\end{array}
\right)\in V$.
For any $y\in{\mathfrak z}_{\mathfrak g}(e)\cap{\cal N}_1$, $y+ke_i\subset{\mathfrak z}_{\mathfrak g}(e_i)\cap{\cal N}_1$; hence there is an action of the additive group ${\mathbb G}_a$ on ${\mathfrak z}_{\mathfrak g}(e_i)\cap{\cal N}_1$ by $\xi\cdot y=y+\xi e_i$.
It follows that each irreducible component of ${\mathfrak z}_{\mathfrak g}(e_i)\cap{\cal N}_1$ is stable under this action of ${\mathbb G}_a$, hence that $x_0+e_i\in V$.
But $\rank(x_0+e_i)=i+\rank(E)>i$ if $E$ is non-zero.
\end{proof}

We therefore consider the subset $Y$ of ${\mathfrak z}_{\mathfrak g}(e_i)\cap{\cal N}_1$ consisting of all $x$ with $E=0$.
The conditions for $x\in{\cal N}_1$ then reduce to: $A^2 =[A,C]+BF = 0$, $AB=0$, $FA=0$.
For each $j$, $0\leq j\leq [i/2]$, let $A_j$ be the $(i\times i)$ matrix $\left(
\begin{array}{ccc}
0 & 0 & I_j \\
0 & 0 & 0 \\
0 & 0 & 0
\end{array}
\right)$, where the top left, top right, bottom left and bottom right submatrices are $j\times j$, the top middle and bottom middle submatrices are $j\times (i-2j)$, the centre left and centre right submatrices are $(i-2j)\times i$, and the central submatrix is $(i-2j)\times (i-2j)$.
Since $Z_G(e_i)$ contains all elements of the form $\left(
\begin{array}{ccc}
g & 0 & 0 \\
0 & 1 & 0 \\
0 & 0 & g
\end{array}
\right)$ with $g\in GL(i)$, it is clear that $Y=\bigcup \overline{Y_j}$, where $$Y_j= Z_G(e_i)\cdot\left\{ \left(
\begin{array}{ccc}
A_j & B & C \\
0 & 0 & F \\
0 & 0 & A_j
\end{array}
\right): A_j B =0, FA_j =0, [A_j,C]+BF = 0\right\}.$$

Moreover, since this is a finite union, each irreducible component of $Y$ is an irreducible component of one of the $\overline{Y_j}$.
Clearly the conditions $A_j B = FA_j = 0$ imply that $B$ and $F$ can be written respectively as $\left(
\begin{array}{c}
B_1 \\
B_2 \\
0
\end{array}
\right)$ and $\left(
\begin{array}{ccc}
0 & F_2 & F_3
\end{array}
\right)$ (where $B_1$ (resp. $B_2$) has $j$ (resp. $(i-2j)$) rows and $F_2$ (resp. $F_3$) has $(i-2j)$ (resp. $j$) columns).
But $A_j\left(\begin{array}{c}
x_1 \\
x_2 \\
x_3
\end{array}\right)=\left(\begin{array}{c}
x_3 \\
0 \\
0
\end{array}\right)$, hence if $x=\left(\begin{array}{c}
0 \\
0 \\
B_1
\end{array}\right)$ then: $$\left(\begin{array}{ccc} I & x & 0 \\ 0 & I & 0 \\ 0 & 0 & I\end{array}\right) \left(
\begin{array}{ccc}
A_j & B & C \\
0 & 0 & F \\
0 & 0 & A_j
\end{array}
\right) \left(\begin{array}{ccc} I & x & 0 \\ 0 & I & 0 \\ 0 & 0 & I\end{array}\right) = \left(
\begin{array}{ccc}
A_j & B+A_j x & C+xF \\
0 & 0 & F \\
0 & 0 & A_j
\end{array}
\right)$$ and $B+A_j x=\left(\begin{array}{c} 0 \\ B_2 \\ 0 \end{array}\right)$.

Similarly, $\left(\begin{array}{lcr}
y_1 & y_2 & y_3
\end{array}\right)A_j=\left(\begin{array}{lcr}
0 & 0 & y_1
\end{array}\right)$.
Hence after a further conjugation we may assume in addition that $F_3=0$.
In other words, any element of $Y_j$ is $Z_G(e_i)$-conjugate to one of the form $\left(\begin{array}{ccc} A_j & B & C \\ 0 & 0 & F \\ 0 & 0 & A_j \end{array}\right)$ such that $B=\left(\begin{array}{c}
0 \\
B_2 \\
0
\end{array}\right)$ and $F=\left(\begin{array}{ccc}
0 & F_2 & 0
\end{array}\right)$.
The equality $[A_j,C]+BF=0$ now implies that $[A_j,C]=0$, $B_2 F_2=0$.
But if $C=\left(\begin{array}{ccc} C_{11} & C_{12} & C_{13} \\ C_{21} & C_{22} & C_{23} \\ C_{31} & C_{32} & C_{33} \end{array}\right)$, where $C_{11}, C_{13},C_{31}$ and $C_{33}$ (resp. $C_{12}$ and $C_{32}$, $C_{21}$ and $C_{23}$, $C_{22}$) are $j\times j$ (resp. $j\times (i-2j)$, $(i-2j)\times j$, $(i-2j)\times (i-2j)$), then $[A_j,C]=\left(\begin{array}{ccc}
C_{31} & C_{32} & C_{11}+C_{33} \\
0 & 0 & C_{21} \\
0 & 0 & C_{31}
\end{array}\right)$.
Hence, conjugating further by $\left(\begin{array}{ccc}
I_i & 0 & z \\
0 & I_{n-2i} & 0 \\
0 & 0 & I_i
\end{array}\right)$, where $z=\left(\begin{array}{ccc} C_{13} & 0 & 0 \\ C_{23} & 0 & 0 \\ C_{11} & C_{12} & 0 \end{array}\right)$, we may assume that $C$ is of the form $\left(
\begin{array}{ccc}
0 & 0 & 0 \\
0 & C_{22} & 0 \\
0 & 0 & 0
\end{array}
\right)$.

We therefore introduce the subset $Y'_j$ of $Y_j$ consisting of all $x$ of the form $\left(
\begin{array}{ccc}
A_j & B & C \\
0 & 0 & F \\
0 & 0 & A_j
\end{array}
\right)$ with $B=\left(
\begin{array}{c}
0 \\
B_2 \\
0
\end{array}
\right)$, $F=\left(
\begin{array}{ccc}
0 & F_2 & 0
\end{array}
\right)$, $C=\left(
\begin{array}{ccc}
0 & 0 & 0 \\
0 & C_{22} & 0 \\
0 & 0 & 0
\end{array}
\right)$, and $B_2 F_2=0$.
We have proved that $Z_G(e_i)\cdot Y_j=Z_G(e_i)\cdot Y_j'$.
Note that the $(n-4j)\times(n-4j)$ matrix $\left(\begin{array}{ccc}
0 & B_2 & C_{22} \\
0 & 0 & F_2 \\
0 & 0 & 0
\end{array}\right)$ is an element of ${\cal N}_1(\mathfrak{gl}(n-4j))$ which commutes with $\left(\begin{array}{ccc} 0 & 0 & I_j \\ 0 & 0 & 0 \\ 0 & 0 & 0 \end{array}\right)$.

\begin{rk}
If we consider pairs $(x,y)\in{\cal C}^{nil}_1$ as modules for $k[X,Y]/(X^2,Y^2)$, then we have shown that for any $y\in Y_j$ the module $M$ corresponding to $(e_i,y)$ can be expressed as $M=W\oplus M'$, where $W$ is a free $k[X,Y]/(X^2,Y^2)$-module of rank $j$ and $M'$ is a submodule of $M$ which is annihilated by $XY$.
\end{rk}

\begin{lemma}\label{harder}
Each irreducible component of $Y$ is properly contained in a closed irreducible subset of ${\mathfrak z}_{\mathfrak g}(e_i)\cap{\cal N}_1$.
\end{lemma}

\begin{proof}
We prove the lemma by induction on $n$.
There is nothing to prove if $n=2$ or $n=3$ (since we assume $n-2i\geq 2$).
By the above remarks, each component of $Y$ is contained in one of the sets $\overline{Y_j}$.

We note that $Y_0={\mathfrak u}\cap{\cal N}_1$, where ${\mathfrak u}$ is the Lie algebra of the unipotent radical of $Z_G(e_i)$.
(Hence $Y_0$ is already closed.)
On the other hand let $j>0$.
Let $e'$ be a nilpotent element of $\mathfrak{gl}(n-4j)$ of partition type $2^{i-2j}.1^{n-2i}$ (in a form as described after Lemma \ref{cpts}), let ${\mathfrak u}'$ be the Lie algebra of the unipotent radical of $Z_{\GL(n-4j)}(e')$ and let $a=\left(\begin{array}{ccc}
A_j & 0 & 0 \\
0 & 0 & 0 \\
0 & 0 & A_j
\end{array}\right)\in{\mathfrak z}_{\mathfrak g}(e_i)$.
We define an injective homomorphism of restricted Lie algebras $\mu:{\mathfrak z}_{\mathfrak{gl}(n-4j)}(e')\rightarrow{\mathfrak z}_{\mathfrak g}(e)$ by $\left(\begin{array}{ccc} A' & B' & C' \\ 0 & E' & F' \\ 0 & 0 & A'\end{array}\right)\mapsto \left(\begin{array}{ccc} A & B & C \\ 0 & E' & F \\ 0 & 0 & A\end{array}\right)$, where $A=\left(\begin{array}{ccc} 0 & 0 & 0 \\ 0 & A' & 0 \\ 0 & 0 & 0 \end{array}\right)$, $C=\left(\begin{array}{ccc} 0 & 0 & 0 \\ 0 & C' & 0 \\ 0 & 0 & 0 \end{array}\right)$, $B=\left(\begin{array}{c} 0 \\ B' \\ 0 \end{array}\right)$ and $F=\left(\begin{array}{ccc} 0 & F' & 0 \end{array}\right)$.
Here the zero submatrices on the top left, top right, bottom left and bottom right (resp. top middle and bottom middle, centre left and centre middle) in $A$ and $C$ are $j\times j$ (resp. $j\times (i-2j)$, $(i-2j)\times j$); those on the top and bottom in $B$ are $j\times (n-2i)$; and those on the left and right in $F$ are $(n-2i)\times j$.
Clearly $a$ commutes with the image of $\mu$ and by the remarks above $a+\mu({\mathfrak u}'\cap{\cal N}_1)=Y'_j$.
But the lemma now follows for $Y_j'$ (and therefore for $\overline{Y_j}=\overline{Z_G(e_i)\cdot Y_j'}$) by the induction hypothesis.
Hence we have only to prove that the statement of the lemma is true for $Y_0$.

Note that $Y_0$ is the set of $x\in{\mathfrak z}_{\mathfrak g}(e_i)$ such that $A=0$, $E=0$ and $BF=0$.
For each $l$ with $0\leq l\leq\min\{ i,n-2i\}$, let $b_l$ be the $i\times (n-2i)$ matrix of the form $\left(
\begin{array}{cc}
I_l & 0 \\
0 & 0
\end{array}
\right)$.
Here the left (resp. right) column is of width $l$ (resp. $n-2i-l$), and the top (resp. bottom) row is of height $l$ (resp. $i-l$).
Since $Z_G(e_i)$ contains all elements of the form $\left(
\begin{array}{ccc}
g & 0 & 0 \\
0 & h & 0 \\
0 & 0 & g
\end{array}
\right)$, it is easy to see that $Y_0=\bigcup_l\overline{{Z_G(e_i)\cdot Z_l}}$, where $$Z_l=\left\{\left(
\begin{array}{ccc}
0 & b_l & C \\
0 & 0 & F \\
0 & 0 & 0
\end{array}
\right)\,:\,b_l F=0\right\}.$$
Moreover, the sets $\overline{Z_G(e_i)\cdot Z_l}$ are clearly irreducible closed subsets of $Y_0$, and $b_l F=0$ if and only if $F$ can be written in the form $\left(\begin{array}{c}
0 \\
f
\end{array}\right)$, where the top part has $l$ rows, and the bottom has $(n-2i-l)$ rows.

Hence it will be enough to prove that, for each $l$, there is a closed irreducible subset of ${\mathfrak z}_{\mathfrak g}(e_i)\cap{\cal N}_1$ which contains $Z_l$ and is not contained in $Y$.
Suppose $0<l<n-2i$, and let $E_0$ be the $(n-2i)\times(n-2i)$ matrix with 1 in the $((n-2i),1)$ position and 0 elsewhere.
Then $b_l E_0=0$ and $E_0\left(\begin{array}{c} 0 \\ f\end{array}\right)=0$, hence the set $\left\{\left(
\begin{array}{ccc}
0 & b_l & C \\
0 & \xi E_0 & F \\
0 & 0 & 0
\end{array}\right)\,:\,\xi\in k,\, b_l F=0\right\}$ is a closed irreducible subset of ${\mathfrak z}_{\mathfrak g}(e_i)\cap{\cal N}_1$ which properly contains $Z_l$.
This proves the lemma in this case.
Let $\theta$ be the automorphism of ${\mathfrak g}$ given by $x\mapsto -J( {^t}x)J^{-1}$, where $J$ is the element of $\GL(n,k)$ with 1 on the antidiagonal, and 0 elsewhere.
Then $\theta(e)=-e$ (hence $\theta$ stabilizes ${\mathfrak z}_{\mathfrak g}(e)\cap{\cal N}_1$) and $\theta$ sends $\overline{Z_G(e)\cdot Z_0}$ into $\overline{Z_G(e)\cdot Z_r}$, where $r=\min\{ n-2i,i\}$.
Hence we have only to prove the statement of the lemma for $Z_{n-2i}$ (assuming therefore that $(n-2i)\leq i$).

Let $b=b_{n-2i}$; we note that left matrix multiplication by $b$ is injective.
Consider the set $U$ of all $x\in{\mathfrak z}_{\mathfrak g}(e_i)\cap{\cal N}_1$ of the form $\left(
\begin{array}{ccc}
A & b & C \\
0 & E & F \\
0 & 0 & A
\end{array}\right)$.
Then the conditions for $x$ to be in ${\cal N}_1$ can be written as: $A^2=[A,C]+bF=0$, $E^2=0$, $Ab=bE$, and $EF=FA$.
But if $Ab=bE$ then $EF=FA$ if and only if $AbF=bFA$.
Since $(\ad A)^2=\ad(A^2)$ it follows that the condition $EF=FA$ is redundant.

If $Ab=bE$ then we can write $A$ in the form $\left(\begin{array}{cc}
E & A_{12} \\
0 & A_{22}
\end{array}\right)$, where $A_{12}\in\Mat_{(n-2i)\times (3i-n)}$ and $A_{22}\in\Mat_{(3i-n)\times (3i-n)}$.
Write $C$ and $F$ respectively in the following forms:
$$C=\left(\begin{array}{cc}
C_{11} & C_{12} \\
C_{21} & C_{22}
\end{array}\right),\;
F=\left(\begin{array}{cc}
F_1 & F_2
\end{array}\right)$$ where the left- (resp. right-) hand columns are of width $(n-2i)$ (resp. $(3i-n)$) and the top (resp. bottom) row of $C$ is of height $(n-2i)$ (resp. $(3i-n)$).
Then $x\in{\cal N}_1$ if and only if $E^2=0$, $A_{22}^2=0$, $EA_{12}=A_{12}A_{22}$ and $$\left(\begin{array}{cc}
F_1 & F_2 \\
0 & 0
\end{array}\right)=\left(\begin{array}{cc}
[E,C_{11}]+A_{12}C_{21} & EC_{12}+C_{12}A_{22}+A_{12}C_{22}+C_{11}A_{12} \\
A_{22}C_{21}+C_{21}E & [A_{22},C_{22}]+C_{21}A_{12}
\end{array}\right)$$
These equalities can be restated as an expression for $F$ in terms of $A,C,E$ together with the conditions: $$A_{22}^2=[A_{22},C_{22}]+C_{21}A_{12}=0,A_{22}C_{21}+C_{21}E=0,E^2=0,EA_{12}+A_{12}A_{22}=0$$
But these conditions are equivalent to: $\left(\begin{array}{ccc}
A_{22} & C_{21} & C_{22} \\
0 & E & A_{12} \\
0 & 0 & A_{22}
\end{array}
\right)^2=0$.
Let $e'$ be a nilpotent element of $\mathfrak{gl}(4i-n)$ of type $2^{3i-n}.1^{n-2i}$, in the form described after Lemma \ref{cpts}.
Then we have proved that there is an isomorphism of affine varieties $U\rightarrow ({\mathfrak z}_{\mathfrak{gl}(4i-n)}(e')\cap{\cal N}_1)\times \Mat_{(n-2i)\times i}$.
Specifically, an element $x$ of the above form is sent to the pair $$\left( \left(\begin{array}{ccc}
A_{22} & C_{21} & C_{22} \\
0 & E & A_{12} \\
0 & 0 & A_{22}
\end{array}
\right) ,\left(\begin{array}{cc}
C_{11} & C_{12}
\end{array}
\right) \right).$$
Notice that $4i-n<n$, and that $(4i-n)-2(3i-n)=n-2i\geq 2$.
Hence it follows by the induction hypothesis that each irreducible component of ${\mathfrak z}_{\mathfrak{gl}(4i-n)}(e')\cap{\cal N}_1$ contains an element such that $E\neq 0$.
Thus the same is true for each irreducible component of $U$.
But $Z_{n-2i}$ is clearly contained in some irreducible component of $U$.
This completes the proof of the lemma.
\end{proof}

We therefore have the required result of this section:

\begin{lemma}\label{lhsreg}
The condition $V_i^{(j)}\subseteq\overline{G\cdot e_i}$ holds if and only if $i=[n/2]$, that is, if and only if $\overline{G\cdot e_i}={\cal N}_1$.
\end{lemma}

\begin{proof}
This follows immediately from Lemmas \ref{Enonzero} and \ref{harder}.
\end{proof}

\section{Equidimensionality}

We proved in the previous section that ${\cal C}^{nil}_1=\overline{G\cdot{\mathfrak z}_{\mathfrak g}(e_m)\cap{\cal N}_1}$, where $m=[n/2]$.
In this section we will show that ${\mathfrak z}_{\mathfrak g}(e_m)\cap{\cal N}_1$ is equidimensional, hence so is ${\cal C}_1^{nil}$.
In fact, we prove equidimensionality of ${\mathfrak z}_{\mathfrak g}(e)\cap{\cal N}_1$ for an arbitrary $e\in{\cal N}_1$.

\begin{lemma}\label{equid2}
Let $e\in{\cal N}_1$ be of partition type $2^i.1^{n-2i}$.

(a) If $n$ is even then ${\mathfrak z}_{\mathfrak g}(e)\cap{\cal N}_1$ has $([i/2]+1)$ irreducible components, each of dimension $(n^2+(n-2i)^2)/4=\dim{\mathfrak z}_{\mathfrak g}(e)/2$.

(b) If $n$ is odd, then ${\mathfrak z}_{\mathfrak g}(e)\cap{\cal N}_1$ has $(i+1)$ irreducible components, each of dimension $(n^2+(n-2i)^2-2)/4=(\dim{\mathfrak z}_{\mathfrak g}(e)-1)/2$.
\end{lemma}

\begin{proof}
As in Sect. 1 we choose $e=\left(\begin{array}{ccc}
0 & 0 & I_i \\
0 & 0 & 0 \\
0 & 0 & 0
\end{array}
\right)$, where the top and bottom rows are of height $i$, the left and right columns are of width $i$, and the central row (resp. middle column) is of height (resp. width) $(n-2i)$.
For $0\leq j\leq [i/2]$ let $A_j$ be the $i\times i$ matrix of the form $\left(\begin{array}{ccc}
0 & 0 & I_j \\
0 & 0 & 0 \\
0 & 0 & 0
\end{array}
\right)$, where the top left, top right, bottom left and bottom right (resp. top middle and bottom middle, centre left and centre right, central) submatrices are $j\times j$ (resp. $j\times (i-2j)$, $(i-2j)\times j$, $(i-2j)\times (i-2j)$).
Similarly, for $0\leq l\leq [(n-2i)/2]$, let $E_l$ be the matrix of the form $\left(\begin{array}{ccc}
0 & 0 & I_l \\
0 & 0 & 0 \\
0 & 0 & 0
\end{array}
\right)$, where the top left, top right, bottom left and bottom right (resp. top middle and bottom middle, centre left and centre right, central) submatrices are $l\times l$ (resp. $l\times (n-2(i+l))$, $(n-2(i+l))\times l$, $(n-2(i+l))\times (n-2(i+l))$).

We denote by $V_{j,l}$ ($0\leq j\leq [i/2],0\leq l\leq[(n-2i)/2]$) the set of all $x\in{\mathfrak z}_{\mathfrak g}(e)\cap{\cal N}_1$ of the form $\left(\begin{array}{ccc}
A_j & B & C \\
0 & E_l & F \\
0 & 0 & A_j
\end{array}
\right)$.
Clearly ${\mathfrak z}_{\mathfrak g}(e)\cap{\cal N}_1=\bigcup_{j,l} \overline{Z_G(e)\cdot V_{j,l}}$, hence each irreducible component of ${\mathfrak z}_{\mathfrak g}(e)\cap{\cal N}_1$ is equal to $\overline{Z_G(e)\cdot V_{j,l}^0}$ for some irreducible component $V_{j,l}^0$ of (some) $V_{j,l}$.
We claim that if $(n-2i)-2l\geq 2$ then no such set $\overline{Z_G(e)\cdot V_{j,l}^0}$ is an irreducible component.

Let $x=\left(\begin{array}{ccc}
A_j & B & C \\
0 & E_l & F \\
0 & 0 & A_j
\end{array}
\right)\in V_{j,l}$, that is, we assume $A_j B=BE_l$, $E_l F=FA_j$ and $[A_j,C]=BF$.
Write $B$ as $\left(\begin{array}{ccc}
B_{11} & B_{12} & B_{13} \\
B_{21} & B_{22} & B_{23} \\
B_{31} & B_{32} & B_{33}
\end{array}
\right)$, where $B_{11}, B_{13},B_{31}$ and $B_{33}$ (resp. $B_{12}$ and $B_{32}$, $B_{21}$ and $B_{23}$, $B_{22}$) are $j\times l$ (resp. $j\times (n-2(i+l))$, $(i-2j)\times l$, $(i-2j)\times (n-2(i+l))$) matrices.
Then $A_j B=\left(\begin{array}{ccc}
B_{31} & B_{32} & B_{33} \\
0 & 0 & 0 \\
0 & 0 & 0
\end{array}
\right)$ and $BE_l= \left(\begin{array}{ccc}
0 & 0 & B_{11} \\
0 & 0 & B_{21} \\
0 & 0 & B_{31}
\end{array}
\right)$.
Hence $A_j B=BE_l$ if and only if $B=\left(\begin{array}{ccc}
B_{11} & B_{12} & B_{13} \\
0 & B_{22} & B_{23} \\
0 & 0 & B_{11}
\end{array}
\right)$.
Similarly, $F=\left(\begin{array}{ccc}
F_{11} & F_{12} & F_{13} \\
0 & F_{22} & F_{23} \\
0 & 0 & F_{11}
\end{array}
\right)$, where $F_{ij}$ has the same dimensions as $B_{ji}^t$.

But $$\left(\begin{array}{ccc} I & y & 0 \\ 0 & I & 0 \\ 0 & 0 & I \end{array}\right) \left(\begin{array}{ccc}
A_j & B & C \\
0 & E_l & F \\
0 & 0 & A_j
\end{array}
\right) \left(\begin{array}{ccc} I & y & 0 \\ 0 & I & 0 \\ 0 & 0 & I \end{array}\right) = \left(\begin{array}{ccc}
A_j & B+A_jy+yE_l & C+yF \\
0 & E_l & F \\
0 & 0 & A_j
\end{array}
\right).$$
Setting $y=\left(\begin{array}{ccc} B_{13} & 0 & 0 \\ B_{23} & 0 & 0 \\ B_{11} & B_{12} & 0 \end{array}\right)$, we see that after conjugating by a suitable element of $Z_G(e)$ we may assume that $B$ is of the form $\left(\begin{array}{ccc}
0 & 0 & 0 \\
0 & B_{22} & 0 \\
0 & 0 & 0
\end{array}
\right)$.
Similarly, conjugating further by a suitable element of the form $\left(\begin{array}{ccc} I & 0 & 0 \\ 0 & I & y \\ 0 & 0 & I \end{array}\right)$ we may assume that $F=\left(\begin{array}{ccc}
0 & 0 & 0 \\
0 & F_{22} & 0 \\
0 & 0 & 0
\end{array}
\right)$.
Under these assumptions it is clear that the condition $[A_j,C]=BF$ implies that $B_{22}F_{22}=0$ and $[A_j,C]=0$.
Finally, conjugating further by an element of the form $\left(\begin{array}{ccc} I & 0 & y \\ 0 & I & 0 \\ 0 & 0 & I \end{array}\right)$, we may assume that $C=\left(\begin{array}{ccc}
0 & 0 & 0 \\
0 & C_{22} & 0 \\
0 & 0 & 0
\end{array}
\right)$, where the zero submatrices on the top left, top right, bottom left and bottom right (resp. top middle and bottom middle, centre left and centre right) are $j\times j$ (resp. $j\times (i-2j)$, $(i-2j)\times j$), and $C_{22}$ is $(i-2j)\times (i-2j)$.
Let $U_{j,l}$ be the set of all $x\in{\mathfrak z}_{\mathfrak g}(e)\cap{\cal N}_1$ with $B,F$, and $C$ of this form.
Then we have proved that $\overline{Z_G(e)\cdot V_{j,l}}=\overline{Z_G(e)\cdot U_{j,l}}$.

Let $a=\left(\begin{array}{ccc} A_j & {} & {} \\ {} & E_l & {} \\ {} & {} & A_j \end{array}\right)\in{\mathfrak{gl}}(n,k)$ and let $e'$ be a nilpotent element of $\mathfrak{gl}(n-4j)$ of type $2^{i-2j}.1^{n-2i}$, in a form analogous to the $e_i$ defined after Lemma \ref{cpts}.
Let ${\mathfrak u}'$ be the Lie algebra of the unipotent radical of $Z_{\GL(n-4j)}(e')$.
There is an injective restricted Lie algebra homomorphism $\mu:{\mathfrak z}_{\mathfrak{gl}(n-4j)}(e')\rightarrow {\mathfrak z}_{\mathfrak g}(e)$ given by $\left(\begin{array}{ccc} A' & B' & C' \\ 0 & E' & F' \\ 0 & 0 & A'\end{array}\right)\mapsto \left(\begin{array}{ccc} A & B & C \\ 0 & E & F \\ 0 & 0 & A\end{array}\right)$, where $A=\left(\begin{array}{ccc} 0 & 0 & 0 \\ 0 & A' & 0 \\ 0 & 0 & 0\end{array}\right)$ such that the zero submatrices on the top left, top right, bottom left and bottom right (resp. top middle and bottom middle, centre left and centre right) are $j\times j$ (resp. $j\times (i-2j)$, $(i-2j)\times j$), and similarly for $B,C,E,F$.
(The dimensions for the corresponding submatrices of $B$ are $j\times l$, $j\times (n-2(i+l))$ and $(i-2j)\times l$; for $C$ are the same as for $A$; for $E$ are $l\times l$, $l\times (n-2(i+l))$ and $(n-2(i+l))\times l$; and for $F$ are $l\times j$, $l\times (i-2j)$ and $(n-2(i+l))\times j$.)

Clearly $a$ commutes with the image of $\mu$, $a+\mu({\mathfrak z}_{\mathfrak{gl}(n-4j)}(e')\cap{\cal N}_1)\subset {\mathfrak z}_{\mathfrak g}(e)\cap{\cal N}_1$ and $a+{\mathfrak u}'\cap{\cal N}_1=U_{j,l}$.
But by Lemma \ref{harder}, ${\mathfrak u}'\cap{\cal N}_1$ is properly contained in a closed irreducible subset of ${\mathfrak z}_{\mathfrak{gl}(n-4j)}(e')\cap{\cal N}_1$.
It follows that $U_{j,l}$ is not an irreducible component of ${\mathfrak z}_{\mathfrak g}(e)\cap{\cal N}_1$ unless $l=[(n-2i)/2]$.

We can now prove equidimensionality.
If $e$ is of type $2^i.1^{2(m-i)}$ then the only sets $\overline{Z_G(e)\cdot V_{j,l}^0}$ which can be irreducible components of ${\mathfrak z}_{\mathfrak g}(e)\cap{\cal N}_1$ are those with $l=(m-i)$.
But now, with the above description of $x\in V_{j,l}$ the possible $B$ and $F$ are:

$$B=\left(\begin{array}{cc}
B_{11} & B_{13} \\
0 & B_{23} \\
0  & B_{11}
\end{array}
\right),\;\; F=\left(\begin{array}{ccc}
F_{11} & F_{12} & F_{13} \\
0 & 0 & F_{11}
\end{array}
\right) .$$

Here $B_{11}$ and $B_{13}$ are $j\times(m-i)$ matrices, $B_{23}$ is $(i-2j)\times (m-i)$, $F_{11}$ and $F_{13}$ are $(m-i)\times j$, and $F_{12}$ is $(m-i)\times (i-2j)$.
The condition $BF=[A_j,C]$ then simply specifies a unique possible $C$ modulo ${\mathfrak z}_{\mathfrak{gl}(i)}(A_j)$.
It follows that $V_{j,(m-i)}$ is irreducible of dimension $2i(m-i)+\dim{\mathfrak z}_{\mathfrak{gl}(i)}(A_j)$.
Let $W_{j,(m-i)}=\overline{Z_G(e)\cdot V_{j,(m-i)}}$ and consider the morphism $\pi:W_{j,(m-i)}\rightarrow \mathfrak{gl}(i)\times\mathfrak{gl}(n-2i)$ given by $\left(\begin{array}{ccc} A & B & C \\ 0 & E & F \\ 0 & 0 & A \end{array}\right)\mapsto (A,E)$.
It is easy to see that the image of $\pi$ is $\overline{\GL(i)\cdot A_j}\times {\cal N}_1(\mathfrak{gl}(n-2i))$.
Moreover, $(\GL(i)\cdot A_j)\times(\GL(n-2i)\cdot E_{m-i})$ is an open subset of $\pi(W_{j,(m-i)})$ and the fibre over each such point is isomorphic to (in fact, conjugate to) $V_{j,(m-i)}$.
It follows by the standard theorem on dimensions (see for example \cite[Thm. 4.3]{hum}) that $\dim W_{j,(m-i)}=i^2-\dim{\mathfrak z}_{\mathfrak{gl}(i)}(A_j)+2(m-i)^2+2i(m-i)+\dim{\mathfrak z}_{\mathfrak{gl}(i)}(A_j)=2m^2-2mi+i^2$.
Since this dimension is independent of $j$, the sets $W_{j,(m-i)}$, $0\leq j\leq [i/2]$ are the irreducible components of ${\mathfrak z}_{\mathfrak g}(e)\cap{\cal N}_1$.
Moreover, $\{ 0\}=\GL(i)\cdot A_0\subset \overline{\GL(i)\cdot A_1}\subset\ldots\subset\overline{\GL(i)\cdot A_{[i/2]}}$, hence $W_{j,(m-i)}$ cannot be contained in $W_{j',(m-i)}$ for $j'<j$.
By equality of dimensions, the $W_{j,(m-i)}$ are distinct.

Similarly, if $e$ is of type $2^i.1^{2(m-i)+1}$ then any irreducible component is of the form $\overline{Z_G(e)\cdot V_{j,(m-i)}^0}$ for some $j$ and some irreducible component $V_{j,(m-i)}^0$ of $V_{j,(m-i)}$.
If $2j<i$ then the possible $B,F$ are of the form:

$$B=\left(\begin{array}{ccc}
B_{11} & b_1 & B_{13} \\
0 & b_2 & B_{23} \\
0 & 0 & B_{11}
\end{array}
\right),\;\; F=\left(\begin{array}{ccc}
F_{11} & F_{12} & F_{13} \\
0 & f_2 & f_3 \\
0 & 0 & F_{11}
\end{array}\right)$$
where $B_{11}$ and $B_{13}$ are $j\times (m-i)$ matrices, $F_{11}$ and $F_{13}$ are $(m-i)\times j$ matrices, $b_1$ (resp. $b_2$) is a column vector of dimension $j$ (resp. $i-2j$), $f_2$ (resp. $f_3$) is a row vector of dimension $i-2j$ (resp. $j$), $B_{23}$ is an $(i-2j)\times (m-i)$ matrix, and $F_{12}$ is $(m-i)\times (i-2j)$.

In this case the condition $BF=[A_j,C]$ specifies a unique value of $C$ modulo ${\mathfrak z}_{\mathfrak{gl}(i)}(A_j)$, and in addition the requirement that $b_2 f_2=0$.
But $b_2$ is a column vector and $f_2$ a row vector, hence $b_2 f_2=0$ implies that either $b_2=0$ or $f_2=0$.
It follows that $V_{j,(m-i)}$ has two irreducible components of equal dimension.
Denote by $V_{j,(m-i)}^+$ the irreducible component defined by $f_2=0$, and by $V_{j,(m-i)}^-$ the irreducible component satisfying $b_2=0$.
Then $\dim V_{j,(m-i)}^+=\dim V_{j,(m-i)}^-=2i(m-i)+i+\dim{\mathfrak z}_{\mathfrak{gl}(i)}(A_j)$.
Let $W_{j,(m-i)}^+=\overline{ Z_G(e)\cdot V_{j,(m-i)}^+}$ and let $W_{j,(m-i)}^-=\overline{Z_G(e)\cdot V_{j,(m-i)}^-}$.
Then by the argument used above, $W_{j,(m-i)}^+$ and $W_{j,(m-i)}^-$ are (irreducible) of dimension $2m^2-2mi+i^2+2m-i=(\dim{\mathfrak z}_{\mathfrak g}(e)-1)/2$.
If $i$ is even and $j=i/2$, then the possible $B,F$ are of the form:

$$B=\begin{pmatrix} B_{11} & b_1 & B_{13} \\ 0 & 0 & B_{11} \end{pmatrix}, F=\begin{pmatrix} F_{11} & F_{13} \\ 0 & f_3 \\ 0 & F_{11}Ê\end{pmatrix}$$

where $B_{11}$ and $B_{13}$ are $j\times (m-i)$ matrices, $F_{11}$ and $F_{13}$ are $(m-i)\times j$ matrices, and $b_1$ (resp. $f_3$) is a column (resp. row) vector of dimension $j$.
Here the condition $BF=[A_j,C]$ merely specifies a unique value of $C$ modulo ${\mathfrak z}_{\mathfrak{gl}(i)}(A_j)$.
It follows that $V_{j,(m-i)}$ is irreducible of dimension $2i(m-i)+i+\dim{\mathfrak z}_{\mathfrak{gl}(i)}(A_j)=i(m-i+1)+\dim{\mathfrak z}_{\mathfrak{gl}(i)}(A_j)$.
Let $W_{i/2,(m-i)}=\overline{Z_G(e)\cdot V_{j,(m-i)}}$.
Then, by exactly the same argument used for the case $j<i/2$, we can see that $\dim W_{i/2,(m-i)}$ is irreducible of dimension $(\dim{\mathfrak z}_{\mathfrak g}(e)-1)/2$.

By equality of dimensions, each irreducible component of ${\mathfrak z}_{\mathfrak g}(e)\cap{\cal N}_1$ is equal to one of the $W_{j,(m-i)}^\pm$ ($0\leq j\leq (i-1)/2$) if $i$ is odd (resp. one of the $W_{j,(m-i)}^\pm$ ($0\leq j\leq i/2-1$) or $W_{i/2,(m-i)}$ if $i$ is even).
Note that there are $(i+1)$ possible choices in either case.
The argument used above for the case where $n$ is even shows that we cannot have $W_{j,(m-i)}^\pm=W_{j',(m-i)}^\pm$ if $j\neq j'$.
Similarly, if $i$ is even then $W_{i/2,(m-i)}\neq W_{j,(m-i)}^\pm$ for $j<i/2$.
Hence it remains to show that $W_{j,(m-i)}^+\neq W_{j,(m-i)}^-$.
By definition $W_{j,(m-i)}^+=\overline{Z_G(e)\cdot V_{j,(m-i)}^+}$.
It is easy to see that $W_{j,(m-i)}^+=W_{j,(m-i)}^-$ is equivalent to: $V_{j,(m-i)}^+\subset Z_G(e)\cdot V_{j,(m-i)}^-$.
Let $x\in V_{j,(m-i)}^+$ and suppose that $g\in Z_G(e)$ satisfies $gxg^{-1}\in V_{j,(m-i)}$.
Then clearly $g$ is of the form $\left(\begin{array}{ccc} h_1 & y_1 & y_2 \\ 0 & h_2 & y_3 \\ 0 & 0 & h_1 \end{array}\right)$ for some $h_1\in Z_{\GL(i)}(A_j)$ and $h_2\in Z_{\GL(2(m-i)+1)}(E_{m-i})$.
Moreover, any $g$ of this form normalizes $V_{j,(m-i)}$.
Let $L$ be the subgroup of $Z_G(e)$ of all elements of this form: $L$ is isomorphic to a product $Z_{\GL(i)}(A_j)\times Z_{\GL(2m-2i+1)}(E_{(m-i)})\times U$, where $U$ is a connected unipotent group.
But therefore $L$ is connected, and hence preserves $V_{j,(m-i)}^+$.
It follows that $W_{j,(m-i)}^+\neq W_{j,(m-i)}^-$.
This completes the proof.
\end{proof}

\begin{rk}
Note that for any irreducible component $V$ of ${\mathfrak z}_{\mathfrak g}(e)\cap{\cal N}_1$ the intersection with the open orbit in ${\cal N}_1$ is non-empty, therefore open.
This is not true for general $p$.
\end{rk}

We now have our result.

\begin{proposition}\label{equid}
Let $k$ be an algebraically closed field of characteristic 2 and let ${\mathfrak g}=\mathfrak{gl}(n,k)$.
Then the restricted nilpotent commuting variety ${\cal C}$ of ${\mathfrak g}$ is equidimensional.

(a) If $n=2m$ then there are $([m/2]+1)$ irreducible components of ${\cal C}$ of dimension $3m^2$.

(b) If $n=2m+1$ then there are $(m+1)$ irreducible components of ${\cal C}$ of dimension $3m(m+1)$.
\end{proposition}

\begin{proof}
By Lemma \ref{lhsreg} each irreducible component of ${\cal C}$ is of the form $\overline{G\cdot(e_m,V)}$ where $V$ is an irreducible component of ${\mathfrak z}_{\mathfrak g}(e_m)\cap{\cal N}_1$.
Let $V$ be such an irreducible component, let $X=\overline{G\cdot(e_m,V)}$ and let $\pi:X\rightarrow{\cal N}_1$ be the restriction to $X$ of the first projection.
Since $G\cdot(e_m,V)$ contains an open subset of its closure, there is an open subset $U$ of ${\cal N}_1$ such that $\dim\pi^{-1}(x)=\dim V\;\forall\, x\in U$.
It follows by the standard theorem on dimensions \cite[Thm. 4.3]{hum} that $\dim X=\dim{\cal N}_1+\dim V=\dim{\mathfrak g}-r$, where $r$ is the codimension of $X$ in ${\mathfrak z}_{\mathfrak g}(e_m)$.
By Lemma \ref{equid2}, ${\mathfrak z}_{\mathfrak g}(e_m)\cap{\cal N}_1$ is equidimensional of dimension $m^2$ (resp. $m^2+m$) if $n=2m$ (resp. $n=2m+1$).
But hence each of the sets $\overline{G\cdot(e_m,V)}$ is an irreducible component of ${\cal C}$, and is of the dimension stated in the proposition.
Moreover, if $X_1=\overline{G\cdot (e_m,V_1)}=\overline{G\cdot(e_m,V_2)}=X_2$ then, since $G\cdot(e_m,V_1)$ contains an open subset of $X_1$, $Z_G(e_m)\cdot V_1$ contains an open subset of $V_2$, and therefore $V_1=V_2$.
This completes the proof of the proposition.

\end{proof}

For later use we now label the components of ${\cal C}$.
Recall from the proof of Lemma \ref{equid2} that if $n=2m$ (resp. $n=2m+1$) then the irreducible components of ${\mathfrak z}_{\mathfrak g}(e_m)\cap{\cal N}_1$ are the sets of the form $\overline{Z_G(e_m)\cdot V_{j,0}}$ (resp. $\overline{Z_G(e_m)\cdot V_{j,0}^\pm}$ and $\overline{Z_G(e_m)\cdot V_{m/2,0}}$ if $m$ is even) for $0\leq j\leq [m/2]$ (resp. $0\leq j<m/2$).
If $n=2m$ then let $X_j=\overline{G\cdot(e_m,V_{j,0})}$.
If $n=2m+1$ then let $X_j^+=\overline{G\cdot(e,V_{j,0}^+)}$, $X_j^-=\overline{G\cdot(e,V_{j,0}^-)}$ and $X_{m/2}=\overline{G\cdot(e,V_{m/2,0})}$ if $m$ is even.

\vspace{8pt}
\begin{rk}\label{notgeneral}
One might ask whether Lemma \ref{lhsreg} or Prop. \ref{equid} is true for fields of arbitrary characteristic.
In fact they both fail.
For example, let $k$ be of characteristic 7 and let $\mathfrak{g}=\mathfrak{gl}(14)$.
Let $e$ be a nilpotent element of ${\mathfrak g}$ of type $7^2$.
Hence $\overline{G\cdot e}={\cal N}_1({\mathfrak g})$.
It is easy to see that ${\mathfrak z}_{\mathfrak g}(e)\cong\mathfrak{gl}(2,k[t]/(t^7))$.
Identify ${\mathfrak z}_{\mathfrak g}(e)$ with $\mathfrak{gl}(2,k[t]/(t^7))$ and write an element $A\in{\mathfrak z}_{\mathfrak g}(e)$ as $A_0+A_1 t+\ldots +A_6 t^6$ where $A_i\in\mathfrak{gl}(2,k)$.
We have $(\sum_0^6 A_i t^i)^7=\sum_0^6 p_i(A) t^i$, where $p_i(A)$ is the sum of all ordered monomials $A_{j_1}A_{j_2}\ldots A_{j_7}$ such that $j_1+j_2+\ldots j_7=i$.
Clearly $x^7=0$ if $A_0=0$.
Hence $t\mathfrak{gl}(2,k[t]/(t^7))$ is a closed irreducible subset of ${\mathfrak z}_{\mathfrak g}(e)\cap{\cal N}_1$ of dimension 24.

On the other hand, if $A_0\neq 0$, then up to conjugacy there is only one possibility such that $A_0^7=0$, namely $A_0=e_{12}$.
Suppose therefore that $A_0=e_{12}$.
We will determine the conditions on the matrices $A_i$ such that $A^7=0$.
We remark first of all that $A_0^2=0$, from which it follows that $p_1(A)=A_0^6A_1+\ldots+A_1A_0^6=0$, and $p_2(A)=A_0^6A_2+\ldots+A_2A_0^6+A_0^5A_1^2+\ldots+A_1^2A_0^5=0$.
The expression for $p_3(A)$ reduces to $(A_0 A_1)^3A_0$.
It follows that $A_1=\left(\begin{array}{cc}
a & b \\
0 & d
\end{array}\right)$ for some $a,b,d\in k$.
Hence $A_0 A_1^j A_0=0$ for any $j\geq 0$.
Inspection of the possible non-zero terms of $p_4$ reveals that $p_4(A)=0$.
Similarly, $p_5(A)$ reduces easily to $A_0^3A_1^3A_2+\ldots+A_2A_1^3A_0^3+A_0^2A_1^5+\ldots+A_1^5A_0^2$.
But each term here either contains $A_0^2$ or it contains $A_0A_1^iA_0$.
Therefore $p_5(A)=0$ also.
Finally, we have $p_6(A)=(A_0A_2)^3A_0+A_0^3A_1^2A_2^2+\ldots+A_2^2A_1^2A_0^3+A_0^2A_1^4A_2+\ldots+A_2A_1^4A_0+(\ad A_1)^6(A_0)$.
Let $A_2=\left(\begin{array}{cc} * & * \\ s & * \end{array}\right)$.
Then $(A_0A_2)^3A_0=s^3e_{12}$, $A_0^3A_1^2A_2^2+\ldots+A_2^2A_1^2A_0^3=-(a-d)^2s^2e_{12}$, $A_0^2A_1^4A_2+\ldots+A_2A_1^4A_0^2=-2(a-d)^4se_{12}$ and $(\ad A_1)^6(A_0)=(a-d)^6e_{12}$.
It follows that $p_6(A)=(s+2(a-d)^2)^3 e_{12}$.
Hence the set $U$ of $A\in{\mathfrak z}_{\mathfrak g}(e)\cap{\cal N}_1$ such that $A_0=e_{12}$ is an irreducible Zariski closed subset of ${\mathfrak z}_{\mathfrak g}(e)$ of dimension 22.
Thus (by the standard theorem on dimensions) $\overline{Z_G(e)\cdot U}$ is an irreducible subset of ${\mathfrak z}_{\mathfrak g}(e)\cap{\cal N}_1$ of dimension 24.
We have proved that ${\mathfrak z}_{\mathfrak g}(e)\cap{\cal N}_1$ has two irreducible components $V_1$ and $V_2$, both of codimension 4 in ${\mathfrak z}_{\mathfrak g}(e)$.
We deduce that the closures $\overline{G\cdot(e,V_1)}$ and $\overline{G\cdot(e,V_2)}$ are irreducible components of ${\cal C}_1^{nil}({\mathfrak g})$  (since $G\cdot e$ is not contained in the closure of any other orbit in ${\cal N}_1$) and are of dimension $(\dim{\mathfrak g}-4)=192$.

On the other hand, let $e'$ be an element of ${\cal N}_1({\mathfrak g})$ of type $7^1.5^1.2^1$, which we may choose to be in Jordan normal form.
Let $T$ be the group of invertible diagonal matrices in $G$ and let $\lambda:k^\times\longrightarrow T$ be the cocharacter such that $e'\in{\mathfrak g}(2;\lambda)$ and the component of $\lambda(t)$ in each Jordan block has determinant 1.
(Hence $\lambda$ is an associated cocharacter for $e'$ in the sense of Pommerening \cite{pom1,pom2}.
In particular, ${\mathfrak z}_{\mathfrak g}(e')\subset\sum_{i\geq 0}{\mathfrak g}(i;\lambda)$.)
Let ${\mathfrak g}(i)={\mathfrak g}(i;\lambda)$ for each $i\in{\mathbb Z}$.
Recall that a {\it toral algebra} ${\mathfrak h}$ is a commutative restricted Lie algebra which has a basis $\{ h_1,\ldots ,h_s\}$ of elements such that $h_i^{[p]}=h_i$. 
It is easily seen that ${\mathfrak z}_{\mathfrak g}(e')\cap{\mathfrak g}(0)$ is a toral algebra of dimension 3 and that ${\mathfrak z}_{\mathfrak g}(e')\cap\sum_{i>0}{\mathfrak g}(i)\subset\sum_{i\geq 2}{\mathfrak g}(i)$.
Hence ${\mathfrak z}_{\mathfrak g}(e')\cap{\cal N}_1={\mathfrak z}_{\mathfrak g}(e')\cap\sum_{i>0}{\mathfrak g}(i)$ is irreducible of codimension 3 in ${\mathfrak z}_{\mathfrak g}(e')$.
It follows that $\overline{G\cdot(e',{\mathfrak z}_{\mathfrak g}(e')\cap{\cal N}_1)}$ is irreducible of dimension $\dim{\mathfrak g}-3=193$.
In particular, it is not contained in either irreducible component of $\overline{G\cdot (e,{\mathfrak z}_{\mathfrak g}(e)\cap{\cal N}_1)}$ (and vice versa, neither irreducible component of $\overline{G\cdot (e,{\mathfrak z}_{\mathfrak g}(e)\cap{\cal N}_1)}$ is contained in $\overline{G\cdot (e',{\mathfrak z}_{\mathfrak g}(e')\cap{\cal N}_1)}$).
Hence in this case neither Lemma \ref{lhsreg} nor equidimensionality of ${\cal C}_1^{nil}$ hold.
\end{rk}

Although Lemma \ref{harder} and Prop. \ref{equid} fail in the above example, we note that the intersection ${\mathfrak z}_{\mathfrak g}(e)\cap{\cal N}_1$ is nevertheless equidimensional.
Indeed, this result appears to be true in general.

\begin{conjecture}
Let $G$ be a reductive group over $k$, and suppose the characteristic of $k$ is good for $G$.
Let $e\in{\cal N}_1$.
Then ${\mathfrak z}_{\mathfrak g}(e)\cap{\cal N}_1$ is equidimensional.
\end{conjecture}

Some laborious but generally straightforward case-checking establishes that the conjecture is true for the cases $G=\GL(4),\GL(5),\SO(5),\GL(6)$.
To illustrate the conjecture and the apparent unpredictability of the number of irreducible components of the intersection ${\mathfrak z}_{\mathfrak g}(e)\cap{\cal N}_1$, we give a few examples.
To determine the irreducible components directly one can use a similar approach to that employed above, that is, consider case-by-case the orbits in ${\mathfrak z}_{\mathfrak g}(e)\cap{\mathfrak g}(0)\cap{\cal N}_1$ where ${\mathfrak g}=\oplus_{i\in{\mathbb Z}}{\mathfrak g}(i)$ is the grading of ${\mathfrak g}$ induced by an associated cocharacter for $e$.

\vspace{8pt}
\noindent{\bf (a) $G=\GL(5,k),\;p=3$.}

(i) $e$ of type $3^1.2^1$.
Here $\dim{\mathfrak z}_{\mathfrak g}(e)=9$.
There are two irreducible components, both of dimension 6.

(ii) $e$ of type $3^1.1^2$.
Then $\dim{\mathfrak z}_{\mathfrak g}(e)=11$, and ${\mathfrak z}_{\mathfrak g}(e)\cap{\cal N}_1$ has two irreducible components, both of dimension 7.

(iii) $e$ of type $2^2.1^1$.
Then $\dim{\mathfrak z}_{\mathfrak g}(e)=13$, and the intersection with ${\cal N}_1$ has three components of dimension 8.

(iv) $e$ of type $2^1.1^3$.
Then $\dim{\mathfrak z}_{\mathfrak g}(e)=17$, and ${\mathfrak z}_{\mathfrak g}(e)\cap{\cal N}_1$ is irreducible of dimension 11.

\vspace{8pt}
\noindent{\bf (b) $G=\GL(6,k),\;p=3$.}

(i) $e$ of type $3^2$.
Then $\dim{\mathfrak z}_{\mathfrak g}(e)=12$ and ${\mathfrak z}_{\mathfrak g}(e)\cap{\cal N}_1$ has two components of dimension 8.

(ii) $e$ of type $3^1.2^1.1^1$.
Then $\dim{\mathfrak z}_{\mathfrak g}(e)=14$ and ${\mathfrak z}_{\mathfrak g}(e)\cap{\cal N}_1$ has five irreducible components of dimension 12.

(iii) $e$ of type $3^1.1^3$.
Then $\dim{\mathfrak z}_{\mathfrak g}(e)=18$ and the intersection with ${\cal N}_1$ is irreducible of dimension 12.

(iv) $e$ of type $2^3$.
Then $\dim{\mathfrak z}_{\mathfrak g}(e)=18$ and ${\mathfrak z}_{\mathfrak g}(e)\cap{\cal N}_1$ has two irreducible components of dimension 12.

For the remaining cases ${\mathfrak z}_{\mathfrak g}(e)\cap{\cal N}_1$ is irreducible.

\vspace{8pt}
\noindent{\bf (c) $G=\Sp(8,k),\;p=3$.}

(i) $e$ of type $3^2.2^2$.
Then $\dim{\mathfrak z}_{\mathfrak g}(e)=12$ and ${\mathfrak z}_{\mathfrak g}(e)\cap{\cal N}_1$ is irreducible of dimension 9.

(ii) $e$ of type $3^2.1^2$.
Then $\dim{\mathfrak z}_{\mathfrak g}(e)=14$ and ${\mathfrak z}_{\mathfrak g}(e)\cap{\cal N}_1$ is irreducible of dimension 10.

(iii) $e$ of type $2^4$.
In this case $\dim{\mathfrak z}_{\mathfrak g}(e)=16$ and ${\mathfrak z}_{\mathfrak g}(e)\cap{\cal N}_1$ is irreducible of dimension 11.

(iv) $e$ of type $2^3.1^2$.
In this case $\dim{\mathfrak z}_{\mathfrak g}(e)=18$ and ${\mathfrak z}_{\mathfrak g}(e)\cap{\cal N}_1$ has two irreducible components, both of dimension 12.

(v) $e$ of type $2^2.1^2$.
Then ${\mathfrak z}_{\mathfrak g}(e)$ is of dimension 22 and ${\mathfrak z}_{\mathfrak g}(e)\cap{\cal N}_1$ has two irreducible components, both of dimension 15.

(vi) $e$ of type $2^1.1^6$.
Then ${\mathfrak z}_{\mathfrak g}(e)$ is of dimension 28 and ${\mathfrak z}_{\mathfrak g}(e)\cap{\cal N}_1$ is irreducible of dimension 19.

\begin{rk}
This conjecture is not true if the characteristic is bad.
Indeed, if $k$ is of characteristic 2, $G$ is simply-connected of type $B_2$ and $e=e_{\alpha_1+\alpha_2}$ (where $\{\alpha_1,\alpha_2\}$ is a basis for the root system of $G$), then ${\mathfrak z}_{\mathfrak g}(e)\cap{\cal N}_1$ has one irreducible component of dimension 3 and one of dimension 4.
On the other hand, the intersection $\Lie(Z_G(e))\cap{\cal N}_1$ is in this case irreducible.
Finally, we remark that the conjecture is not true for arbitrary nilpotent elements.
For example, if $p=2$ and $e$ is an element of $\mathfrak{gl}(6)$ of type $3^2$, then ${\mathfrak z}_{\mathfrak g}(e)\cap{\cal N}_1$ has two components, one of dimension 8 and one of dimension 6.
\end{rk}

\section{Indecomposable components}

We remarked in the introduction that the variety ${\cal C}={\cal C}_1^{nil}$ can be considered as the variety of $n$-dimensional modules for the group algebra $k\Gamma$, where $\Gamma$ is a product of two cyclic groups of order 2.
More generally, let $A$ be any finitely generated associative algebra.
To give an $r$-dimensional vector space $V$ the structure of an $A$-module is simply to give a homomorphism $A\rightarrow\End(V)$.
Such a homomorphism is determined by the values on a set of generators for $A$.
Hence, on choosing a basis for $V$, the set $\Mod_A^r(k)$ of possible $A$-module structures on $k^r$ embeds as a Zariski closed subset of the product of a finite number of copies of $\Mat_r(k)$.
The general linear group $\GL(r,k)$ acts on $\Mod_A^r(k)$ by simultaneous conjugation on the coordinates.
It is clear that two points of $\Mod_A^r(k)$ are $\GL(r,k)$-conjugate if and only if the corresponding modules are isomorphic.

An irreducible component of $\Mod_A^r(k)$ is called {\it indecomposable} if all points in an open subset correspond to indecomposable modules.
A version of the Krull-Remak-Schmidt Theorem holds for irreducible components of $\Mod_A^r(k)$.
Let $C_i:1\leq i\leq l$ be irreducible components of varieties of $A$-modules $\Mod_A^{r_i}(k)$ and assume that $r_1+r_2+\ldots r_l=r$.
There is a morphism $\GL(r,k)\times C_1\times C_2\times\ldots\times C_l\rightarrow\Mod_A^r(k)$.
Denote the closure of the image by $\overline{C_1\oplus C_2\oplus\ldots\oplus C_l}$.
Then any irreducible component of $\Mod_A^r(k)$ can be expressed in an essentially unique way as a direct sum $\overline{C_1\oplus C_2\oplus\ldots\oplus C_l}$ of indecomposable components $C_i$ of module varieties $\Mod_A^{r_i}(k)$ (originally proved in \cite{pena}; see also \cite{cbs}. Such a direct sum is not always an irreducible component for arbitrary $C_i$: see \cite[Thm. 1.2]{cbs}.)
Here we express each component of ${\cal C}$ as a direct sum of indecomposable components.

Recall that the irreducible components of ${\cal C}$ are labelled $X_j$ if $n=2m$ with $0\leq j\leq [m/2]$ (resp. $X_j^\pm$ ($0\leq j<m/2$) if $n=2m+1$ and $m$ is odd, $X_j^\pm$ ($0\leq j<m/2$) and $X_{m/2}$ if $m$ is even).
For arbitrary $r$, we denote by $X_j(\mathfrak{gl}(r))$ or $X_j^\pm(\mathfrak{gl}(r))$ the irreducible components of ${\cal C}_1^{nil}(\mathfrak{gl}(r))$ described in this way.
Let $W$ be the irreducible component of ${\cal C}_1^{nil}(\mathfrak{gl}(4))$ given by $W=\overline{\GL(4)\cdot(e_{12}+e_{34},e_{13}+e_{24})}$ (the ``free component of rank 1"), and let $U$ be the irreducible component of ${\cal C}_1^{nil}(\mathfrak{gl}(2))$ given by $U=\{ (ae,be)\,:\, e^{[2]}=0,\; a,b\in k\}$.
It is easy to see that $W$ (resp. $U$) is an indecomposable component of ${\cal C}_1^{nil}(\mathfrak{gl}(4))$ (resp. ${\cal C}_1^{nil}(\mathfrak{gl}(2))$).
If $A=k\Gamma$ for some group $\Gamma$ and $M$ is any left $A$-module, then the dual vector space $M^*$ has the structure of a left-$A$-module with $g\in\Gamma$ acting via $(g\cdot\chi)(m)=\chi(g^{-1}\cdot m)$ for each $\chi\in M^*$, $m\in M$.
In these circumstances, if $V$ is an irreducible component of $\Mod_A^r$, then we denote by $V^*$ the dual component $\{ M^*:M\in V\}$.
Let $\triv={\cal C}_1^{nil}(\mathfrak{gl}(1))$ denote the variety of one-dimensional $A$-modules.
Clearly $\triv$ clearly consists of a single point.
We have the following:

\begin{proposition}

(a) Suppose $n=2m$.
Then $X_j\cong \overline{W^j\oplus U^{n-4j}}$ and $X_j^*=X_j$.

(b) Suppose $n=2m+1$.
Then $X_0^\pm$ are indecomposable components of ${\cal C}$.
Moreover, $X_j^+=\overline{W^j\oplus X_0^+(\mathfrak{gl}(n-4j))}$, $X_j^-=\overline{W^j\oplus X_0^-(\mathfrak{gl}(n-4j))}$ and $(X_j^+)^*=X_j^-$.
If $m$ is even, then $X_{m/2}=\overline{W^{m/2}\oplus \triv}$.
Moreover, $X_{m/2}^*=X_{m/2}$.
\end{proposition}

\begin{proof}
Let $e=e_m$ and let $\theta$ be the automorphism of ${\mathfrak g}$ given by $x\mapsto -J(^{t}x)J^{-1}$, where $J$ is the matrix with 1 on the anti-diagonal, and 0 elsewhere.
Let $G$ (resp. $\theta$) act diagonally on ${\mathfrak g}\times{\mathfrak g}$, hence on each irreducible component of ${\cal C}$.
Since $\theta$ is a restricted Lie algebra automorphism of ${\mathfrak g}$, its induced action on ${\cal C}$ permutes the irreducible components of ${\cal C}$.
Clearly $\theta(X_j)=X_j^*$ (resp. $\theta(X^\pm_j)=(X^\pm_j)^*$).
If $n$ is even, recall that $X_j=\overline{G\cdot(e,V_{j,0})}$, where $V_{j,0}$ is the set defined in the proof of Lemma \ref{equid2}.
Since the zero here is superfluous, we will write $V_j$ for $V_{j,0}$.
Similarly, we will write $V_j^\pm$ for $V_{j,0}^\pm$ in the case $n$ odd below and $V_{m/2}$ for $V_{m/2,0}$ if $m$ is even.
Clearly $\theta(e)=-e$ and $\theta(\diag(A_j,A_j))=-\diag(A_j,A_j)$.
It is easy to choose $g$ such that $\Ad g(e)=-e$ and $\Ad g(\diag(A_j,A_j))=-\diag(A_j,A_j)$.
Hence $\Ad g\circ\theta(X_j)=X_j$, and thus $X_j^*=X_j$.
This proves the second statement of (a).
We proved in Lemma \ref{equid2} that $\overline{Z_G(e)\cdot V_j}=\overline{Z_G(e)\cdot U_j}$, where $U_j$ is the set of elements of ${\mathfrak z}_{\mathfrak g}(e)$ of the form $\left(\begin{array}{cc}
A_j & B \\
0 & A_j
\end{array}\right)$, $B=\left(\begin{array}{ccc}
0 & 0 & 0 \\
0 & B_{22} & 0 \\
0 & 0 & 0
\end{array}\right)$.
Here the zero submatrices on the top left, top right, bottom left and bottom right (resp. top middle and bottom middle, centre left and centre right) of $B$ are $j\times j$ (resp. $j\times (m-2j)$, $(m-2j)\times j$) and $B_{22}$ is $(m-2j)\times (m-2j)$.
It follows at once that $X_j\cong \overline{W^j\oplus X_0(\mathfrak{gl}(n-4j))}$.
Hence we have only to prove (a) for the case $j=0$.
Thus consider $V_0=\left\{\left(\begin{array}{cc}
0 & B \\
0 & 0
\end{array}\right):B\in\Mat_{m\times m}(k)\right\}$.
The set of semisimple elements is dense in ${\mathfrak{gl}}(m)$, hence the subset of $V_0$ of elements such that $B$ is semisimple is dense.
But any such element is $Z_G(e)$-conjugate to one such that $B$ is diagonal.
Hence $X_0$ is the closure of $G\cdot \left\{ \left(\left(\begin{array}{cc} 0 & I \\ 0 & 0 \end{array}\right),\left(\begin{array}{cc} 0 & B \\Ê0 & 0 \end{array}\right) \right):\; B\mbox{ diagonal}\right\}$.
This proves (a).

For (b), suppose first of all that $m$ is even.
Recall from the proof of Lemma \ref{equid2} that $X_{m/2}=\overline{Z_G(e)\cdot U_{m/2}}$, where $U_{m/2}$ is the set of $x\in{\mathfrak z}_{\mathfrak g}(e)$ of the form $\begin{pmatrix}ÊA_j & 0 & 0 \\ 0 & 0 & 0 \\ 0 & 0 & A_j\end{pmatrix}$.
Here the top and bottom (resp. middle) rows are of height $m$ (resp. 1) and the left and right (resp. central) columns are of width $m$ (resp. 1).
But then clearly $X_{m/2}\cong \overline{W^m\oplus\triv}$.
Since $W=W^*$ and clearly $\triv=\triv^*$, we have also that $X_{m/2}\cong X_{m/2}^*$.

Consider therefore the components $X_j^\pm$ (for arbitrary $m$).
It is easy to see that $\theta(e)=-e$ and $\theta\left(\begin{array}{ccc} A_j & b & C \\ 0 & 0 & f \\ 0 & 0 & A_j\end{array}\right)=\left(\begin{array}{ccc} -A_j & -J(^tf) & -J ^t CJ \\ 0 & 0 & -^t bJ \\ 0 & 0 & -A_j\end{array}\right)$.
We recall from the proof of Lemma \ref{equid2} that $\overline{Z_G(e)\cdot V_j^+}=\overline{Z_G(e)\cdot U_j^+}$, where $U_j^+$ is the set of all $x\in{\mathfrak z}_{\mathfrak g}(e)$ of the form
$\left(\begin{array}{ccc}
A_j & b & C \\
0 & 0 & 0 \\
0 & 0 & A_j
\end{array}\right)$: $b$ is a column vector such that the first $j$ and the last $j$ entries are zero, and $C$ is of the form $\begin{pmatrix}Ê0 & 0 & 0 \\ 0 & C_{22} & 0 \\ 0 & 0 & 0 \end{pmatrix}$, where zero submatrices at the top left, top right, bottom left and bottom right (resp. top middle and bottom middle, centre left and centre right) are of dimension $j\times j$ (resp. $j\times (i-2j)$, $(i-2j)\times j$) and $C_{22}$ is $(i-2j)\times (i-2j)$.
Similarly, $\overline{Z_G(e)\cdot V_j^-}=\overline{Z_G(e)\cdot U_j^-}$, where $U_j^-$ is the set of all $x\in{\mathfrak z}_{\mathfrak g}(e)$ of the form $\begin{pmatrix}ÊA_j & 0 & C \\Ê0 & 0 & f \\ 0 & 0 & A_j \end{pmatrix}$.
Here $f$ is a row vector of dimension $i$ such that the first and last $j$ entries are zero.
Hence, after applying conjugation by a suitably chosen element $g$, $\Ad g\circ\theta(V_j^+)=V_j^-$.
It follows that $(X_j^+)^*=X_j^-$ and $(X_j^-)^*=X_j^+$.
As above, the argument in the proof of Lemma \ref{equid2} shows that $X_j^+=\overline{W^j\oplus X_0^+(\mathfrak{gl}(n-4j))}$, and similarly for $X_j^-$.
Hence we have only to prove that $X_0^+$ is indecomposable (since the result for $X_0^-$ follows on taking the dual).

Let $x=e_{1,m+1}+e_{2,m+2}+\ldots+e_{m,2m}$.
Clearly $x\in V_0^+$.
Moreover, $Z_G(x)\cap Z_G(e)$ is the set of all elements of the form $\left(\begin{array}{ccc} aI_m & y & z \\ 0 & a & 0 \\ 0 & 0 & aI_m\end{array}\right)$.
In particular, $\dim Z_G(x)\cap Z_G(e)=m^2+m+1$.
It follows that $\dim \overline{G\cdot (e,x)}=3(m^2+m)=\dim X_0^+$, hence $\overline{G\cdot (e,x)}=X_0^+$.
To show that $X_0^+$ is indecomposable, it will therefore suffice to show that the module corresponding to $(e,x)$ is indecomposable.
But $Z_G(e)\cap Z_G(x)$ contains no non-trivial idempotents.
This completes the proof.
\end{proof}

\end{document}